\documentclass[12pt]{article}

\usepackage{amssymb}
\usepackage{amsmath}
\usepackage{cite}
\usepackage{dsfont}
\usepackage[arrow, matrix, curve]{xy}

\begin{document}

\renewcommand{\citeleft}{{\rm [}}
\renewcommand{\citeright}{{\rm ]}}
\renewcommand{\citepunct}{{\rm,\ }}
\renewcommand{\citemid}{{\rm,\ }}

\newcounter{abschnitt}
\newtheorem{satz}{Theorem}
\newtheorem{coro}[satz]{Corollary}
\newtheorem{theorem}{Theorem}[abschnitt]
\newtheorem{koro}[theorem]{Corollary}
\newtheorem{prop}[theorem]{Proposition}
\newtheorem{lem}[theorem]{Lemma}
\newtheorem{expls}[theorem]{Examples}
\newtheorem{rem}[theorem]{Remark}
\setcounter{abschnitt}{0}

\renewenvironment{quote}{\list{}{\leftmargin=0.62in\rightmargin=0.62in}\item[]}{\endlist}

\newcounter{saveeqn}
\newcommand{\alpheqn}{\setcounter{saveeqn}{\value{abschnitt}}
\renewcommand{\theequation}{\mbox{\arabic{saveeqn}.\arabic{equation}}}}
\newcommand{\reseteqn}{\setcounter{equation}{0}
\renewcommand{\theequation}{\arabic{equation}}}

\hyphenpenalty=9000

\sloppy

\phantom{a}

\vspace{-2.5cm}

\begin{center}
\begin{Large} {\bf Sharp Sobolev Inequalities via Projection Averages} \\[0.6cm] \end{Large}

\begin{large} Philipp Kniefacz and Franz E.\ Schuster\end{large}
\end{center}

\vspace{-0.9cm}

\begin{quote}
\footnotesize{ \vskip 1cm \noindent {\bf Abstract.}	
A family of sharp $L^p$ Sobolev inequalities is established by averaging the length of $i$-dimensional projections of the gradient of a function.
Moreover, it is shown that each of these new inequalities directly implies the classical $L^p$~Sobolev inequality of Aubin and Talenti and that the strongest member of this family is the only affine invariant one among them -- the affine $L^p$~Sobolev inequality of Lutwak, Yang, and Zhang. When $p = 1$, the entire family of new Sobolev inequalities is extended to functions of bounded variation to also allow for a complete classification of all extremal functions in this case.
}
\end{quote}

\vspace{0.6cm}

\centerline{\large{\bf{ \setcounter{abschnitt}{1}
\arabic{abschnitt}. Introduction}}}

\alpheqn

\vspace{0.6cm}

The fruitful interplay between analysis and geometry is probably highlighted most prominently by the rich theory of Sobolev inequalities and, in particular, by
its best known representative -- the \emph{sharp} $L^p$ Sobolev inequality in $\mathbb{R}^n$. While the latter is often stated for functions from the Sobolev space $W^{1,p}(\mathbb{R}^n)$ (consisting of $L^p$ functions with weak $L^p$ partial derivatives), a more natural setting for it is the larger \emph{homogeneous Sobolev space} (see, e.g., \textbf{\cite[\textnormal{Chapter 11}]{leoni}}) defined by
\[\dot{W}^{1,p}(\mathbb{R}^n) := \{f \in L^{p^*}(\mathbb{R}^n): \nabla f \in L^p(\mathbb{R}^n)\},  \]
where $1 \leq p < n$ and $p^*=np/(n-p)$. The \emph{$L^p$ Sobolev inequality} states that if $f \in \dot{W}^{1,p}(\mathbb{R}^n)$, then
\begin{equation} \label{lpsob}
\left ( \int_{\mathbb{R}^n}\|\nabla f(x)\|^p\,dx \right )^{1/p} \geq a_{n,p}\, \|f\|_{p^*},
\end{equation}
where $\|\cdot\|$ denotes the standard Euclidean norm on $\mathbb{R}^n$ and we write $\|f\|_p$ for the usual $L^p$ norm of $f$ in $\mathbb{R}^n$.
The exact value of the optimal constant $a_{n,p}$ (see below) was first computed for $p > 1$ by Aubin \textbf{\cite{aubin1976problemes}} and, independently, by Talenti \textbf{\cite{talenti1976best}}, who made critical use of the classical isoperimetric inequality to reduce (\ref{lpsob}) to a $1$-dimensional problem. In the case $p = 1$, it was previously shown by Maz'ya \textbf{\cite{mazya}} and Federer and Fleming \textbf{\cite{fedflem}} that the sharp $L^1$ Sobolev inequality is actually equivalent to the isoperimetric inequality.

While the explicit knowledge of the optimal constant has proven beneficial in certain areas of mathematical physics, its importance is far outweighed by the classification of the extremal functions in (\ref{lpsob}). Apparently, it was known for some time that with the help of a rearrangement inequality of Brothers and Ziemer \textbf{\cite{brothziem88}} all extremizers could be identified. However, the first explicit and selfcontained proof that equality holds in (\ref{lpsob}) for $p > 1$ if and only if there exist $a, b>0$, and $x_0\in\mathbb{R}^n$ such that
\begin{equation} \label{equallpsob}
f(x)=\pm\left (a+b\|(x-x_0)\|^{p/(p-1)}\right )^{1-n/p}
\end{equation}
was given by Cordero-Erausquin et al.\ \textbf{\cite{cordnazvill}} (and in a more general form). They also pointed out the disadvantage of considering inequality (\ref{lpsob}) merely for functions in $W^{1,p}(\mathbb{R}^n)$, since its extremizers do not belong to that space when $p \geq \sqrt{n}$. A similar problem arises in the case $p = 1$, where it was known for some time that the natural setting for (\ref{lpsob}) is the space $BV(\mathbb{R}^n)$ of functions of bounded variation, since here the extremizers are characteristic functions of Euclidean balls (see below). For more information on the history of the $L^p$ Sobolev inequality and its ongoing prominent role in different areas, we refer to \textbf{\cite{brothziem88, cianchietal2009, cordnazvill, fuscoetal, leoni}} and the references therein.

\vspace{0.2cm}

It was a major breakthrough when in 2002, Lutwak, Yang, and Zhang \textbf{\cite{lutwak2002sharp}} (building on earlier work of Zhang \textbf{\cite{zhang1999affine}}) established the first \emph{affine invariant} $L^p$~Sobolev inequality. They discovered that replacing the length of the gradient in (\ref{lpsob}) by an average of a suitable power of the length of $1$-dimensional projections of the gradient leads to a significantly stronger inequality than (\ref{lpsob}). More precisely, it was shown in \textbf{\cite{lutwak2002sharp}} that if $1 \leq p < n$ and
$f \in \dot{W}^{1,p}(\mathbb{R}^n)$, then
\begin{equation}\label{lpaffSob}
\left ( \int_{\mathrm{Gr}_{n,1}} \left ( \int_{\mathbb{R}^n} \|\nabla f(x)|E\|^p\, dx \right )^{-n/p} dE \right )^{-1/n} \geq c_{n,p}\, \|f\|_{p^*},
\end{equation}
where we denote by $\mathrm{Gr}_{n,j}$ the Grassmannian of $j$-dimensional subspaces of $\mathbb{R}^n$ and write $\nabla f(x)|E$ for the orthogonal projection of $\nabla f(x)$ to $E \in \mathrm{Gr}_{n,j}$. Throughout, integration over $\mathrm{Gr}_{n,j}$ is with respect to the invariant probability measure on $\mathrm{Gr}_{n,j}$.
The optimal constant $c_{n,p}$ in (\ref{lpaffSob}) was explicitly determined in \textbf{\cite{zhang1999affine}} for $p = 1$ and in \textbf{\cite{lutwak2002sharp}} for $p > 1$.
Later, Wang \textbf{\cite{wang2013affine}} and, independently, Nguyen \textbf{\cite{nguyen}} proved that equality holds in (\ref{lpaffSob}) for $p > 1$ if and only if
\begin{equation} \label{equ1742}
f(x) = \pm \left (a + \|A(x - x_0)\|^{p/(p-1)} \right )^{1-n/p},
\end{equation}
for some $a > 0$, $A \in \mathrm{GL}(n)$, and $x_0 \in \mathbb{R}^n$. For $p = 1$, Wang \textbf{\cite{wang2012affine}} extended (\ref{lpaffSob}) to functions from $BV(\mathbb{R}^n)$ and showed that the extremals in this case are precisely the characteristic functions of ellipsoids. The proof of (\ref{lpaffSob}) by Lutwak~et~al.\ as well as the characterization of extremals by Wang rely on an affine isoperimetric inequality, known as the $L^p$~Petty projection inequality, established in \textbf{\cite{petty1971isoperimetric}} for $p = 1$ and in \textbf{\cite{lutwak2000lp}} for $p > 1$ (see Section 2). Nguyen used a different but equivalent affine inequality, known as the $L^p$~Busemann--Petty centroid inequality, going back to \textbf{\cite{Petty:1961}} for $p = 1$ and \textbf{\cite{campi2002lp, lutwak2000lp}} for $p > 1$, that was first used in this context by Haddad et al.\ \textbf{\cite{haddad2016}}.

\vspace{0.2cm}

The impact of the affine $L^p$~Sobolev inequality (\ref{lpaffSob}) is virtually unparalleled in convex geometric analysis, as it constitutes the seminal result in a rapidly evolving theory of affine analytic inequalities (see, e.g., \textbf{\cite{chiachietal, haberl2009asymmetric, haberl2012asymmetric, haddad2018, haddad2019, Lin, ludwigxiazhang}}). Among this theory's most recent achievements is a large family of sharp $L^p$~Sobolev inequalities by Haberl and the second author \textbf{\cite{haberl2018affine}} that had not just the classical inequality (\ref{lpsob}) and the affine $L^p$~Sobolev inequality (\ref{lpaffSob}) as special cases but also an $(n-1)$-dimensional counterpart to (\ref{lpaffSob}): If $1 \leq p < n$ and $f \in \dot{W}^{1,p}(\mathbb{R}^n)$, then
\begin{equation}\label{lphssob}
\left ( \int_{\mathrm{Gr}_{n,n-1}} \left ( \int_{\mathbb{R}^n} \|\nabla f(x)|E\|^p\, dx \right )^{-n/p} dE \right )^{-1/n} \geq \widetilde{c}_{n,p}\, \|f\|_{p^*}
\end{equation}
with equality for $p > 1$ if and only if $f(x)$ is of the form (\ref{equallpsob}). In the case $p = 1$, inequality (\ref{lphssob}) was extended to $BV(\mathbb{R}^n)$ in \textbf{\cite{haberl2018affine}} and it was shown that equality holds precisely for characteristic functions of Euclidean balls (see below, for the value of the optimal constant). While (\ref{lphssob}) is \emph{not} affine invariant, it was proved in \textbf{\cite{haberl2018affine}} that it is stronger and directly implies the classical $L^p$ Sobolev inequality (\ref{lpsob}) in the same way as the affine $L^p$ Sobolev inequality (\ref{lpaffSob}) does. However, among these three inequalities the affine invariant one was shown in \textbf{\cite{haberl2018affine}} to be the strongest one.

\vspace{0.2cm}

Comparing inequalities (\ref{lpsob}), (\ref{lpaffSob}), and (\ref{lphssob}) raises two natural questions:
\begin{itemize}
\item Is there a family of (sharp) $L^p$ Sobolev inequalities obtained by averaging the length of \emph{$i$-dimensional} gradient projections that unifies (\ref{lpsob}), (\ref{lpaffSob}), and (\ref{lphssob})?
\item How are these gradient projection Sobolev inequalities related to each other?
\end{itemize}

The purpose of this article is to answer both these questions. In order to state our results it is convenient to introduce the following notation.

\vspace{0.3cm}

\noindent {\bf Definition} \emph{Suppose that $1 \leq i \leq n$ and $1 \leq p < n$. For $f \in \dot{W}^{1,p}(\mathbb{R}^n)$, we define
\begin{equation} \label{defeip}
\mathcal{E}_{i,p}(f) = \left ( \int_{\mathrm{Gr}_{n,i}} \left ( \frac{2\omega_{i+p-2}}{i\omega_i \omega_{p-1}} \int_{\mathbb{R}^n} \|\nabla f(x)|E\|^p\,dx \right )^{-n/p}dE\right )^{-1/n},
\end{equation}
where $\omega_p = \pi^{p/2}/\Gamma\left (1 + \frac{p}{2}\right)$.}

\vspace{0.3cm}

With our first main result we give an affirmative answer to the first of the above questions.

\begin{satz} \label{satz:sobolev_lp}
Suppose that $1 \leq i \leq n$ and $1 \leq p < n$. If $f\in \dot{W}^{1,p}(\mathbb{R}^n)$, then
\begin{equation} \label{eq:sobolev_lp}
\mathcal{E}_{i,p}(f)  \geq c_{n,p} \|f\|_{p^*},
\end{equation}
where
\begin{equation*}
c_{n,p} = \left( \frac{2\omega_{n+p-2}}{\omega_{n}\omega_{p-1}} \right)^{1/p}\left(\frac{n-p}{p-1}\right)^{1-1/p}
\left( \frac{\omega_n\Gamma( \frac np )\Gamma(n+1-\frac np)}{\Gamma(n)} \right)^{1/n}.
\end{equation*}
For $p > 1$, equality holds if and only if $f(x)$ has the form (\ref{equallpsob}) when $i > 1$, and if and only if $f(x)$ has the form (\ref{equ1742}) when $i = 1$.
\end{satz}

In Section 4, we actually prove a Sobolev inequality more general than (\ref{eq:sobolev_lp}), where the Euclidean norm of the gradient projection in (\ref{defeip}) can be replaced by any norm whose unit ball is a polar zonoid in $E$ (see Section 2 for definitions).

\pagebreak

The answer to the second question is provided by our next theorem, which shows that the functionals $\mathcal{E}_{i,p}$ form a decreasing sequence.

\begin{satz} 	\label{satz:sobolev_chain_lp}
Suppose that $1 \leq i \leq n$ and $1 \leq p < n$. If $f\in \dot{W}^{1,p}(\mathbb{R}^n)$, then
\begin{equation*}
\mathcal{E}_{n,p}(f) \geq \mathcal{E}_{n-1,p}(f) \geq \cdots \geq \mathcal{E}_{2,p}(f) \geq \mathcal{E}_{1,p}(f).
\end{equation*}
\end{satz}

Note that for any $i \leq n - 1$, inequality (\ref{eq:sobolev_lp}) directly implies the classical $L^p$~Sobolev inequality by Theorem \ref{satz:sobolev_chain_lp}, and that the affine $L^p$ Sobolev inequality (\ref{lpaffSob}) is stronger than any of the inequalities from (\ref{eq:sobolev_lp}) for $i \geq 2$.

\vspace{0.2cm}

While for $p = 1$, the Sobolev inequalities (\ref{eq:sobolev_lp}) are still sharp, their extremizers no longer belong to the space $\dot{W}^{1,1}(\mathbb{R}^n)$ but rather to the space of functions of bounded variation. In this setting the classical Sobolev inequality states that if $f \in BV(\mathbb{R}^n)$, then
\begin{equation} \label{L1Sob}
\|Df\| \geq n\omega_n^{1/n} \|f\|_{\frac{n}{n-1}},
\end{equation}
where the vector valued Radon measure $Df$ is the weak gradient of $f$ and $\|Df\|$ denotes its total variation in $\mathbb{R}^n$ (see Section 2). Equality holds in (\ref{L1Sob}) if and only if $f$ is a multiple of the characteristic function of a Euclidean ball.

In order to extend Theorems \ref{satz:sobolev_lp} and \ref{satz:sobolev_chain_lp} for $p = 1$ to functions of bounded variation, we introduce the following notation.

\vspace{0.25cm}

\noindent {\bf Definition} \emph{Suppose that $1 \leq i \leq n$. For $f \in BV(\mathbb{R}^n)$, we define
\begin{equation*}
	\mathcal{E}_i(f) = \left (\int_{\mathrm{Gr}_{n,i}}\left ( \frac{2\omega_{i-1}}{i\omega_i} \int_{\mathbb{R}^n}\| \sigma_f | E \| \,d|Df| \right ) ^{-n}dE\right)^{-1/n},
\end{equation*}
where $|Df|$ denotes the variation measure of $Df$ and $\sigma_f$ the Radon--Nikodym derivative of $Df$ with respect to $|Df|$.}

\vspace{0.25cm}

The aforementioned extensions can now be conveniently stated as follows.

\begin{satz} \label{satz:sobolev}
Suppose that $1\leq i \leq n$. If $f \in BV(\mathbb{R}^n)$, then
\begin{equation} \label{bvchainsob}
\mathcal{E}_i(f) \geq \frac{2\omega_{n-1}}{\omega_n^{1-1/n} }\, \|f\|_{\frac{n}{n-1}}
\end{equation}
with equality if and only if $f$ is the multiple of a characteristic function of a ball when $i > 1$ and that of an ellipsoid when $i = 1$. Moreover,
\begin{equation} \label{bvstrongest}
\mathcal{E}_n(f) \geq \mathcal{E}_{n-1}(f) \geq \cdots \geq \mathcal{E}_2(f) \geq \mathcal{E}_1(f).
\end{equation}
\end{satz}

In the case $i = 1$, inequality (\ref{bvchainsob}) reduces to the affine invariant Zhang--Sobolev inequality from \textbf{\cite{wang2012affine, zhang1999affine}} which, by (\ref{bvstrongest}), is the strongest of the Sobolev inequalities provided by Theorem \ref{satz:sobolev}. In particular, it directly implies the case $i = n - 1$ of (\ref{bvchainsob}) obtained in \textbf{\cite{haberl2018affine}}, which, in turn, is stronger than the classical Sobolev inequality for functions of bounded variation (\ref{L1Sob}) -- the case $i = n$ of (\ref{bvchainsob}).

\pagebreak

\centerline{\large{\bf{ \setcounter{abschnitt}{2}
\arabic{abschnitt}. Background material}}}

\reseteqn \alpheqn \setcounter{theorem}{0} \vspace{0.6cm}

\newcommand{\M}{\mathcal{M}}
\newcommand{\C}{\mathcal{C}}

In this section we recall basic notions and results from the $L^p$ Brunn--Minkowski theory of convex bodies as well as some definitions and facts about functions of bounded variation required in the proofs of our main results. For additional details on the material presented here, we refer to the excellent monograph \textbf{\cite{schneider2013convex}} by Schneider and the classic text \textbf{\cite{evansmeasure}} by Evans and Gariepy.

\vspace{0.2cm}

The setting for this article is $n$-dimensional Euclidean space $\mathbb{R}^n$, where we always assume that $n \geq 3$.
We denote by $\mathcal{K}^n$ the space of convex bodies (that is, non-empty, compact, convex sets) in $\mathbb{R}^n$ with the Hausdorff metric.
If $K \in \mathcal{K}^n$ has non-empty interior, then we write $|K|$ for its volume. The \emph{polar set} of a closed, convex subset $K \subseteq \mathbb{R}^n$ containing the origin is defined by $K^{\circ} = \{x \in \mathbb{R}^n: x \cdot y \leq 1 \mbox{ for all } y \in K\}$. If $K \in \mathcal{K}^n$ contains the origin in its interior, then $K^{\circ\circ} = K$ and $K^{\circ} \in \mathcal{K}^n$ is called the \emph{polar body} of $K$. Next recall that each $K \in \mathcal{K}^n$ is uniquely determined by the values of its \emph{support function} $h(K,x) = \max\{x \cdot y: y \in K\}$ for $x \in \mathbb{R}^n$. Clearly,
\begin{equation} 	\label{eq:support_function_thetainv}
h(\vartheta K,x) = h(K,\vartheta^{-1}x), \qquad x \in \mathbb{R}^n,
\end{equation}
for every $\vartheta \in \mathrm{SO}(n)$. For an origin-symmetric, closed, and convex $K \subseteq \mathbb{R}^n$, let $\|x\|_K = \min \{\lambda \geq 0: x \in \lambda K\}$, $x \in \mathbb{R}^n$, denote the \emph{Minkowski functional} of $K$ and note that if $K \in \mathcal{K}^n$ has non-empty interior, then $\|\cdot\|_K$ is the norm with unit ball $K$. Moreover, for every origin-symmetric $K \in \mathcal{K}^n$,
\begin{equation*} 
h(K,\cdot) = \|\cdot\|_{K^{\circ}}.
\end{equation*}

Next let $1 \leq i \leq n - 1$ and $E, F \in \mathrm{Gr}_{n,i}$ be given and choose $\vartheta \in \mathrm{SO}(n)$ such that $F = \vartheta E$. If $K \in \mathcal{K}^n$ is origin-symmetric and $i$-dimensional such that $K \subseteq E$, we write $K(F)$ instead of $\vartheta K$ for the rotated copy of $K$ contained in $F$. In this case, it is easy to see that for every $x \in \mathbb{R}^n$,
\begin{equation} \label{suppprojection}
\|x|F\|_{K(F)^{\circ}} = h(\vartheta K,x).
\end{equation}

While the classical Brunn--Minkowski theory of convex bodies emerges from combining the notion of volume with that of Minkowski addition,
the development of its more modern $L^p$ extension, initiated by Lutwak \textbf{\cite{lutwak1993brunn, lutwak1996brunn}}, is a result of merging volume with
the $L^p$ Minkowski addition of convex bodies. To make this more explicit, let $1 \leq p < \infty$ and suppose that $K, L \in \mathcal{K}^n$ contain the origin in their interiors.
For $t > 0$, the \emph{$L^p$ Minkowski combination} $K +_p t\cdot L \in \mathcal{K}^n$ was defined in \textbf{\cite{firey1962p}} by
\[h(K +_p t\cdot L,\cdot)^p = h(K,\cdot)^p + t\,h(L,\cdot)^p.  \]
Note that when $p = 1$, we have $K +_1 t\cdot L = K + tL = \{x + ty: x \in K, y \in L\}$.

It was shown in \textbf{\cite{lutwak1993brunn}} by Lutwak that to each $K \in \mathcal{K}^n$ containing the origin in its interior one can associate a unique Borel measure $S_p(K,\cdot)$ on $\mathbb{S}^{n-1}$, the \emph{$L^p$ surface area measure} of $K$, such that
\[\lim_{t\rightarrow 0^+} \frac{|K +_p t\cdot L| - |K|}{t} = \frac{1}{p} \int_{\mathbb{S}^{n-1}} h(L,u)^p\,dS_p(K,u)  \]
for every $L \in \mathcal{K}^n$ containing the origin in its interior. Moreover, $S_p(K,\cdot)$ is absolutely continuous with respect to the classical surface area measure
$S_1(K,\cdot) = S(K,\cdot)$ and its Radon--Nikodym derivative is $h(K,\cdot)^{1-p}$.

The \emph{Cauchy projection formula} relates the projection function of a convex body $K \in \mathcal{K}^n$ with the cosine transform of its surface area measure in the following way,
\begin{equation} \label{cauchyproj}
\mathrm{vol}_{n-1}(K|u^\perp) = \frac 12\int_{\mathbb{S}^{n-1}} |u \cdot v|\, dS(K, v), \qquad u \in \mathbb{S}^{n-1}.
\end{equation}
Noting that the right-hand side of (\ref{cauchyproj}) is sublinear in $u$ and therefore a support function, Minkowski used it to define the \emph{projection body} of $K$ by
\begin{equation} \label{defl1projbod}
h(\Pi K, x) = \frac 12\int_{\mathbb{S}^{n-1}} |x \cdot v|\, dS(K, v), \qquad x \in \mathbb{R}^n.
\end{equation}

Based on (\ref{defl1projbod}), a natural $L^p$ extension of the projection body operator was defined by Lutwak, Yang, and Zhang in \textbf{\cite{lutwak2000lp}}. For $1\leq p< \infty$ and $K \in \mathcal{K}^n$ containing the origin its interior, the \emph{$L^p$ projection body} of $K$ is given by
\begin{equation} \label{deflpprojbod}
h(\Pi_p K,x)^p = \frac{\omega_{p-1}}{2\omega_{n+p-2}}\int_{\mathbb{S}^{n-1}} |x\cdot v|^p \,dS_p(K,v), \qquad x \in \mathbb{R}^n,
\end{equation}
where the normalizing constant was chosen such that $\Pi_p \mathbb{B}^n = \mathbb{B}^n$ for the Euclidean unit ball $\mathbb{B}^n$ in $\mathbb{R}^n$. Note that when $p = 1$, (\ref{deflpprojbod}) is well defined for all $K \in \mathcal{K}^n$ and that, in this case, $\Pi_1 K = \omega_{n-1}^{-1}\Pi K$.

The range of the $L^p$ projection body map is contained in the class of \emph{$L^p$~zonoids}. For $p \geq 1$, an origin-symmetric convex body $K \in \mathcal{K}^n$ is an $L^p$~zonoid if and only if there exists an even measure $\mu$ on $\mathbb{S}^{n-1}$ (which is uniquely determined when $p$ is not an even integer) such that
\begin{equation} \label{zonoidmeas}
h(K,x)^p = \int_{\mathbb{S}^{n-1}} |x \cdot v|^p\,d\mu(v), \qquad x \in \mathbb{R}^n.
\end{equation}
In the following, we denote by $Z_p^{\mu}$ the $L^p$~zonoid generated in this way by $\mu$. $L^1$~zonoids are usually just called zonoids and we simply write $Z^{\mu}$ instead of $Z_1^{\mu}$ for the zonoid generated by $\mu$. For more information on $L^p$ zonoids and the various contexts, where they arise naturally, we refer to \textbf{\cite[\textnormal{Chapter~3.5}]{schneider2013convex}}.

\pagebreak

The fundamental affine isoperimetric inequality for $L^p$ projection bodies was established by Lutwak, Yang, and Zhang \textbf{\cite{lutwak2000lp}} and is known as the \emph{$L^p$ Petty projection inequality} (see also \textbf{\cite{boeroezky2013, campi2002lp, haberl2009general, LYZ2010a}} for alternative proofs). It states that if $1 \leq p < \infty$ and $K \in \mathcal{K}^n$ contains the origin in its interior, then
\begin{equation} 	\label{eq:lp_ppi}
|\Pi_p^{\circ} K||K|^{\frac{n-p}p} \leq \omega_n^{n/p}
\end{equation}
with equality for $p > 1$ if and only if $K$ is an ellipsoid centered at the origin. Here and in the following, we write $\Pi_p^{\circ}K$ instead of $(\Pi_p K)^\circ$.
The case $p = 1$ of (\ref{eq:lp_ppi}) was already obtained in 1971 by Petty \textbf{\cite{petty1971isoperimetric}} who proved the inequality for all $K \in \mathcal{K}^n$ with non-empty interior and showed that equality holds in this case if and only if $K$ is an ellipsoid (for recent extensions to non-convex sets, see \textbf{\cite{wang2012affine, zhang1999affine}}).

\vspace{0.2cm}

As already mentioned in the introduction, the $L^p$ Petty projection inequality (\ref{eq:lp_ppi}) was crucial in the first proofs of the affine $L^p$ Sobolev inequality (\ref{lpaffSob}). Conversely, (\ref{lpaffSob}) is a functional form of (\ref{eq:lp_ppi}) in the sense that the choice of a suitable function $f$ in the Sobolev inequality (\ref{lpaffSob}) allows to recover the isoperimetric inequality (\ref{eq:lp_ppi}) (see (\ref{PiLYZbod}) and Proposition \ref{prop:lyz_results} below). The correspondence between these inequalities was made even more evident by a new more conceptualized proof of the affine $L^p$~Sobolev inequality by Lutwak, Yang, and Zhang \textbf{\cite{lutwak2006optimal}}, where they associate to each $f \in \dot{W}^{1,p}(\mathbb{R}^n)$ a convex body
$\langle f \rangle_p$. This convexification of a Sobolev function is the content of the following theorem. (Note that measurable functions on $\mathbb{R}^n$ that coincide almost everywhere with respect to Lebesgue measure are considered equal.)

\begin{theorem} \label{the:lyz_volume_normalized} \emph{(\!\!\textbf{\cite{lutwak2006optimal}})} If $1 \leq p < \infty$ and $f \in \dot{W}^{1,p}(\mathbb{R}^n)$ is not identically $0$, then there exists a unique origin-symmetric convex body $\langle f \rangle_p$ with non-empty interior such that
\begin{equation*}
\int_{\mathbb{R}^n} g(\nabla f(x))^p\,dx = \frac{1}{\left | \langle f \rangle_p \right |} \int_{\mathbb{S}^{n-1}} g(u)^p\,dS_p\left ( \langle f \rangle_p,u \right )
\end{equation*}
for every even continuous function $g: \mathbb{R}^n \rightarrow [0,\infty)$ that is positively $1$-homogeneous.
\end{theorem}

In order to see how to apply Theorem \ref{the:lyz_volume_normalized} in our context, note that if $f \in \dot{W}^{1,p}(\mathbb{R}^n)$ and we define
$K = \left |\langle f \rangle_p\right |^{-1/(n-p)} \langle f \rangle_p$, then, by (\ref{deflpprojbod}) and Theorem \ref{the:lyz_volume_normalized},
\begin{equation} \label{PiLYZbod}
h(\Pi_pK,y)^p = \frac{\omega_{p-1}}{2\omega_{n+p-2}} \int_{\mathbb{R}^n} |\nabla f(x) \cdot y|^p \,dx, \qquad y \in \mathbb{R}^n.
\end{equation}
Hence, by the polar coordinate formula for volume, the left-hand side of (\ref{lpaffSob}) coincides up to a constant with $|\Pi_p^\circ K|^{-1/n}$. Consequently, the $L^p$ Petty projection inequality reduces the proof of the affine $L^p$ Sobolev inequality (\ref{lpaffSob}) to a sharp estimate of the volume $|\langle f \rangle_p|$ in terms of $\|f\|_{p^*}$ (which was established in \textbf{\cite{lutwak2006optimal}}). For our purposes, this viewpoint will be helpful to settle the equality cases in Theorem \ref{satz:sobolev_lp}. In order to establish the equality conditions in
Theorem \ref{satz:sobolev}, we require an extension of Theorem~\ref{the:lyz_volume_normalized} when $p = 1$ to functions of bounded variations obtained by Wang \textbf{\cite{wang2012affine}}. To this end let us first recall a few basic facts about the space $BV(\mathbb{R}^n)$.

\pagebreak

A function $f \in L^1(\mathbb{R}^n)$ belongs to $BV(\mathbb{R}^n)$ if for every $1\leq i\leq n$, there exists a finite signed Radon measure $D_if$ on $\mathbb{R}^n$ such that
\begin{equation} \label{bvradonmeas}
\int_{\mathbb{R}^n} f \frac{\partial\xi}{\partial x_i}\,dx = -\int_{\mathbb{R}^n} \xi \,dD_if
\end{equation}
for all compactly supported $C^1$ functions $\xi$ on $\mathbb{R}^n$. A subset $L \subseteq \mathbb{R}^n$ is called a \emph{set of finite perimeter} if $\mathds{1}_L \in BV(\mathbb{R}^n)$.

The \textit{variation} $|Df|$ of the vector valued measure $Df=(D_1f,\ldots,D_nf)$ on $\mathbb{R}^n$ is the non-negative Radon measure whose value at a Borel set $L \subseteq \mathbb{R}^n$ is given by
\[|Df|(L) = \sup_{\pi} \sum_{A \in \pi} |Df(A)|,  \]
where the supremum is taken over all partitions $\pi$ of $L$ into a countable number of disjoint measurable subsets. For $f \in BV(\mathbb{R}^n)$, let $\sigma_f$ denote the Radon-Nikodym derivative of $Df$ with respect to $|Df|$. Then, by (\ref{bvradonmeas}),
\begin{equation*}
\int_{\mathbb{R}^n} f\,\mathrm{div}\phi\,dx = - \int_{\mathbb{R}^n} \phi\cdot \sigma_f\,d|Df|
\end{equation*}
for  all continuously differentiable vector fields $\phi$ on $\mathbb{R}^n$ with compact support.

\vspace{0.2cm}

The aforementioned extension of Theorem~\ref{the:lyz_volume_normalized} to functions of bounded variation by Wang~\textbf{\cite{wang2012affine}} can now be stated as follows.

\begin{theorem} \label{the:lyz_bounded_variation} \emph{(\!\! \textbf{\cite{wang2012affine}})} If $f \in BV(\mathbb{R}^n)$ is not identically $0$, then there exists a unique origin-symmetric convex body $\langle f \rangle$ with non-empty interior such that
\[\int_{\mathbb{R}^n} g(\sigma_f(x))\,d|Df|(x) = \int_{\mathbb{S}^{n-1}} g(u)\,dS\left ( \langle f \rangle,u \right )  \]
for every even continuous function $g: \mathbb{R}^n \rightarrow \mathbb{R}$ that is positively $1$-homogeneous.
\end{theorem}

In case the level sets of $f \in \dot{W}^{1,p}(\mathbb{R}^n)$ or $BV(\mathbb{R}^n)$ are all homothets of a fixed convex body, this body's shape is recovered by $\langle f \rangle_p$ or $\langle f \rangle$, respectively. To make this more precise, let $\mu_f: [0,\infty) \rightarrow [0,\infty]$, $\mu_f(t) = \left | \{x \in \mathbb{R}^n: |f(x)|> t\} \right |$, denote the \emph{distribution function} of $f$ and recall that for an origin-symmetric $K \in \mathcal{K}^n$ with non-empty interior, the \emph{convex symmetrization} of $f$ with respect to $K$ is the function $f^K: \mathbb{R}^n \rightarrow [0,\infty]$ given by
\[f^K(x) = \inf \left \{t > 0: \mu_f(t) \leq \omega_n\|x\|_{\widetilde{K}}^n \right \},  \]
where $\widetilde{K}$ is the dilate of $K$ such that $|\widetilde{K}| = \omega_n$.

\begin{prop} \label{prop:lyz_results} \emph{(\!\!\textbf{\cite{wang2012affine, wang2013affine}})}
Let $K \in \mathcal{K}^n$ be an origin-symmetric convex body with non-empty interior.
\begin{enumerate}
\item[(a)] $\langle \mathds{1}_K \rangle = K$.
\item[(b)] If $1 \leq p < \infty$ and $f \in \dot{W}^{1,p}(\mathbb{R}^n)$, then $\langle f^K \rangle_p$ is a dilate of $K$.
\end{enumerate}
\end{prop}

\pagebreak

\vspace{1cm}
\centerline{\large{\bf{ \setcounter{abschnitt}{3}
\arabic{abschnitt}. Auxiliary results}}}

\reseteqn \alpheqn \setcounter{theorem}{0}
\vspace{0.6cm}

Here, we first recall how to lift integration of functions and measures on the homogeneous spaces $\mathbb{S}^{n-1}$ and $\mathrm{Gr}_{n,i}$ to the Lie group $\mathrm{SO}(n)$ and use this in the second part to prove the underlying geometric inequality behind Theorem \ref{satz:sobolev_chain_lp} and (\ref{bvstrongest}).

\vspace{0.2cm}

From now on, let $\{e_1, \ldots, e_n\}$ denote a fixed orthonormal basis of $\mathbb{R}^n$ and for $1 \leq i \leq n$, let $E_i \in \mathrm{Gr}_{n,i}$ and $S^{i-1} \subseteq \mathbb{S}^{n-1}$ be given by
\[E_i = \mathrm{span}\{e_1, \ldots, e_i\} \qquad \mbox{and} \qquad S^{i-1} = \mathbb{S}^{n-1} \cap E_i.  \]
We write $\mathrm{SO}(i)$ for the subgroup of $\mathrm{SO}(n)$ which leaves $E_i$ invariant and acts as the identity on $E_i^{\perp}$. Note that for $2 \leq i \leq n$, $\mathrm{SO}(i)$ acts transitively on $S^{i-1}$ and that
\[\mathrm{SO}(1) \subseteq \mathrm{SO}(2) \subseteq \cdots \subseteq \mathrm{SO}(n - 1) \subseteq \mathrm{SO}(n).  \]

Next, recall that the unit sphere $\mathbb{S}^{n-1}$ is a homogeneous space with respect to the action of $\mathrm{SO}(n)$. Therefore, $\mathbb{S}^{n-1}$ is diffeomorphic to $\mathrm{SO}(n)/\mathrm{SO}(n-1)$ and there is a one-to-one correspondence between functions and measures on $\mathbb{S}^{n-1}$ and right $\mathrm{SO}(n-1)$ invariant measures on $\mathrm{SO}(n)$. More precisely, if $\mu$ is a measure on $\mathbb{S}^{n-1}$, then there exists a unique right $\mathrm{SO}(n-1)$ invariant measure $\breve{\mu}$ on $\mathrm{SO}(n)$ such that
\begin{equation} \label{mulift}
\int_{\mathbb{S}^{n-1}} f(u)\,d\mu(u) = \int_{\mathrm{SO}(n)} f(\phi e_n)\,d\breve{\mu}(\phi)
\end{equation}
for every $f \in C(\mathbb{S}^{n-1})$. In other words, the pushforward of $\breve{\mu}$ under the natural projection $\pi: \mathrm{SO}(n) \rightarrow \mathbb{S}^{n-1}$, $\pi(\phi)=\phi e_n$, is $\mu$ (see, e.g., \textbf{\cite{grinberg1999convolutions, haberl2018affine}} for more details).

Since, similarly, for $2 \leq i \leq n$, $S^{i-1}$ is diffeomorphic to $\mathrm{SO}(i)/\mathrm{SO}(i-1)$, any measure on $\mathbb{S}^{n-1}$ whose support is concentrated on $S^{i-1}$ may be lifted either to a right $\mathrm{SO}(n-1)$ invariant measure on $\mathrm{SO}(n)$ or to a right $\mathrm{SO}(i-1)$ invariant measure on $\mathrm{SO}(i)$. In particular, we make frequent use of the fact that if $\sigma_i$ denotes the restriction of the $(i-1)$-dimensional Hausdorff measure to $S^{i-1}$, then
\begin{equation} \label{sigmaisoi}
\int_{\mathbb{S}^{n-1}} f(u)\,d\sigma_i(u) = i\omega_i \int_{\mathrm{SO}(i)} f(\phi e_i)\,d\phi
\end{equation}
for every $f \in C(\mathbb{S}^{n-1})$, where integration on the right is with respect to the Haar probability measure on $\mathrm{SO}(i)$.

Since the Lie group $\mathrm{SO}(n)$ also acts transitively on $\mathrm{Gr}_{n,i}$ for every $1 \leq i \leq n - 1$, the Grassmannian $\mathrm{Gr}_{n,i}$ is diffeomorphic to $\mathrm{SO}(n)/\mathrm{S}(\mathrm{O}(i) \times \mathrm{O}(n-i))$, where the subgroup $\mathrm{S}(\mathrm{O}(i) \times \mathrm{O}(n-i))$ is the stabilizer of $E_i$ in $\mathrm{SO}(n)$. Thus, we can also lift integration with respect to measures on $\mathrm{Gr}_{n,i}$ to the group $\mathrm{SO}(n)$. Specifically, we have (as for the sphere $\mathbb{S}^{n-1}$) that for every $f \in C(\mathrm{Gr}_{n,i})$,
\begin{equation} \label{mulifttograss}
\int_{\mathrm{Gr}_{n,i}} f(E)\,dE = \int_{\mathrm{SO}(n)} f(\phi E_i)\,d\phi.
\end{equation}

\pagebreak

Let us turn to $L^p$ zonoids. For $p \geq 1$ and an even measure $\mu$ on $\mathbb{S}^{n-1}$, the support function of the $L^p$ zonoid $Z_p^{\mu}$ generated by $\mu$ can be
written, by (\ref{zonoidmeas}) and (\ref{mulift}), as
\begin{equation} \label{eq:rotation_one_dim}
h(Z_p^{\mu}, x)^p = \int_{\mathrm{SO}(n)}\!\! |x\cdot\phi e_n|^p\,d\breve\mu(\phi), \qquad x \in \mathbb{R}^n.
\end{equation}
Since for $x \in \mathbb{R}^n$, we have
\[\|x\|^p = h(\mathbb{B}^n,x)^p =  \frac{\omega_{p-1}}{2\omega_{n+p-2}}\int_{\mathbb{S}^{n-1}}\! |x\cdot u|^p\,d\sigma_n(u) = \frac{n\omega_n \omega_{p-1}}{2\omega_{n+p-2}}
\int_{\mathrm{SO}(n)}\!\! |x\cdot\phi e_n|^p\,d\phi, \]
we see that the Euclidean unit ball $\mathbb{B}^n$ is an $L^p$ zonoid for any $p \geq 1$. More general, if we denote by $D^i = \mathbb{B}^n \cap E_i$ the $i$-dimensional unit ball in $E_i$,
then, by (\ref{sigmaisoi}),
\begin{equation} \label{normdisclpzon}
\|x|E_i\|^p = h(D^i,x)^p =  \frac{\omega_{p-1}}{2\omega_{i+p-2}}\int_{\mathbb{S}^{n-1}}\! |x\cdot u|^p\,d\sigma_i(u) = \frac{i\omega_i \omega_{p-1}}{2\omega_{i+p-2}}
\int_{\mathrm{SO}(i)}\!\! |x\cdot\phi e_i|^p\,d\phi,
\end{equation}
which shows that also $D^i$ is an $L^p$ zonoid for every $2 \leq i \leq n$ and any $p \geq 1$. When $i = 1$, we have $h(D^1,x)^p = h([-e_1,e_1],x)^p = |x \cdot e_1|^p = \|x|E_1\|^p$, that is, $D^1$ is also and $L^p$ zonoid for any $p \geq 1$.

For $1 \leq i \leq n$ and $p \geq 1$, we now define
\[q_{i,p} = \frac{2\omega_{i+p-2}}{i\omega_i \omega_{p-1}} \qquad \mbox{and} \qquad \nu_{i,p} = i\omega_iq_{i,p}\,\sigma_i.  \]
Then the measure $\nu_{i,p}$ on $\mathbb{S}^{n-1}$ is concentrated on $S^{i-1}$ and, by (\ref{normdisclpzon}), the $L^p$ zonoid~$D^i_p$ generated by $\nu_{i,p}$ satisfies for every $x \in \mathbb{R}^n$,
\begin{equation} \label{yoda1717}
h(D_p^i,x)^p = \int_{\mathbb{S}^{n-1}}\! |x\cdot u|^p\,d\nu_{i,p}(u) = \int_{\mathrm{SO}(i)}\!\! |x\cdot\phi e_i|^p\,d\phi  = q_{i,p}\,\|x|E_i\|^p.
\end{equation}

Noting that $D^1_p = [-e_1,e_1]$ for every $p \geq 1$ and $h(D^1_p,x) = |x \cdot e_1| = \|x|E_1\|$, we conclude from the invariance of the Haar measure on $\mathrm{SO}(i)$ and (\ref{yoda1717}) that for any $2 \leq i \leq n$ and every $x \in \mathbb{R}^n$,
\begin{equation} \label{eq:rotation_support_function}
q_{i,p}\,\|x|E_i\|^p = \int_{\mathrm{SO}(i)}\!\! |x\cdot\phi e_i|^p\,d\phi = \int_{\mathrm{SO}(i)}\!\! |x\cdot\phi e_1|^p\,d\phi = \int_{\mathrm{SO}(i)}\!\! \|x|\phi E_1\|^p\,d\phi.
\end{equation}

With the next lemma, we generalize (\ref{eq:rotation_support_function}) by proving a useful relation between the length of $j$-dimensional projections in terms of averages of the length of $i$-dimensional projections, when $i < j$.

\begin{lem} \label{lem:rotational_norm} If $p \geq 1$ and $1\leq i < j \leq n$, then
	\begin{equation*}
		q_{j,p}\, \|x|E_j\|^p = q_{i,p}\int_{\mathrm{SO}(j)} \| x |\phi E_i\|^p\,d\phi
	\end{equation*}
for every $x \in \mathbb{R}^n$.	
\end{lem}
\noindent {\it Proof.} First note that, by (\ref{yoda1717}), the desired relation is equivalent to
\begin{equation} \label{equivstat1742}
h(D^j_p, x)^p = \int_{\mathrm{SO}(j)}\!\! h(\phi D^i_p, x)^p\,d\phi.
\end{equation}
In order to prove (\ref{equivstat1742}), we use (\ref{eq:support_function_thetainv}) and a combination of (\ref{yoda1717}) and (\ref{eq:rotation_support_function}) to see that
\[\int_{\mathrm{SO}(j)}\!\! h(\phi D^i_p, x)^p\,d\phi = \int_{\mathrm{SO}(j)}\!\! h(D^i_p,\phi^{-1} x)^p\,d\phi =
\int_{\mathrm{SO}(j)} \int_{\mathrm{SO}(i)}\!\! h(\theta D^1_p,\phi^{-1} x)^p\,d\theta    \,d\phi. \]
Thus, from an application of Fubini's theorem, (\ref{eq:support_function_thetainv}), and the fact that $\mathrm{SO}(i) \subseteq \mathrm{SO}(j)$ as well as the invariance of the Haar measure on $\mathrm{SO}(j)$, we obtain
\[\int_{\mathrm{SO}(j)}\!\! h(\phi D^i_p, x)^p\,d\phi = \int_{\mathrm{SO}(i)} \int_{\mathrm{SO}(j)}\!\! h(\phi\theta D^1_p,x)^p\,d\phi\,d\theta = \int_{\mathrm{SO}(j)}\!\! h(\phi D^1_p,x)^p\,d\phi.  \]
Finally, another application of (\ref{yoda1717}) and (\ref{eq:rotation_support_function}) completes the proof of (\ref{equivstat1742}). \hfill $\blacksquare$

\vspace{0.3cm}

With the help of Lemma \ref{lem:rotational_norm}, we can now prove a geometric inequality which is critical for our proofs of Theorem \ref{satz:sobolev_chain_lp} and (\ref{bvstrongest}).

\begin{theorem} \label{the:geometric_inequality}
Suppose that $1 \leq i < j \leq n$ and $p \geq 1$. If $K \in \mathcal{K}^n$ contains the origin in its interior, then
\[ \int_{\mathrm{Gr}_{n,j}}\!\!\! \left (\! q_{j,p}\!\int_{\mathbb{S}^{n-1}}\!\!\! \|u|E\|^p \, dS_p(K,u)\! \right)^{\!-n/p}\!dE \leq
\int_{\mathrm{Gr}_{n,i}}\!\!\! \left (\! q_{i,p}\!\int_{\mathbb{S}^{n-1}}\!\!\! \|u|F\|^p \, dS_p(K,u)\! \right)^{\!-n/p}\!dF.  \]
\end{theorem}

\noindent {\it Proof.} By Lemma \ref{lem:rotational_norm} and Fubini's Theorem, we have that, for any fixed $\theta \in \mathrm{SO}(n)$,
\[\left (\! q_{j,p}\!\int_{\mathbb{S}^{n-1}}\!\!\! \|u|\theta E_j\|^p \, dS_p(K,u)\! \right)^{\!-n/p} =
\left (\! q_{i,p}\!\int_{\mathrm{SO}(j)}\int_{\mathbb{S}^{n-1}}\!\!\! \|u|\theta\phi E_i\|^p \, dS_p(K,u)\,d\phi\! \right)^{\!-n/p}.\]
Hence, by Jensen's inequality,
\[\left (\! q_{j,p}\!\int_{\mathbb{S}^{n-1}}\!\!\! \|u|\theta E_j\|^p \, dS_p(K,u)\! \right)^{\!-n/p} \leq
\int_{\mathrm{SO}(j)}\!\! \left (\! q_{i,p}\!\int_{\mathbb{S}^{n-1}}\!\!\! \|u|\theta\phi E_i\|^p \, dS_p(K,u)\! \right)^{\!-n/p}\! d\phi.\]
Integrating now both sides of this inequality with respect to the Haar probability measure on $\mathrm{SO}(n)$, followed by an application of Fubini's theorem and the invariance of the Haar measure on the right-hand side, we obtain
\[\int_{\mathrm{SO}(n)}\!\! \left (\! q_{j,p}\!\int_{\mathbb{S}^{n-1}}\!\!\! \|u|\theta E_j\|^p \, dS_p(K,u)\! \right)^{\!-n/p}\!d\theta \leq
\int_{\mathrm{SO}(n)}\!\! \left (\! q_{i,p}\!\int_{\mathbb{S}^{n-1}}\!\!\! \|u|\theta E_i\|^p \, dS_p(K,u)\! \right)^{\!-n/p}\!\! \,d\theta  \]
which concludes the proof by (\ref{mulifttograss}). \hfill $\blacksquare$

\pagebreak

\centerline{\large{\bf{\setcounter{abschnitt}{4}
\arabic{abschnitt}. Proofs of the main results}}}

\reseteqn\alpheqn\setcounter{theorem}{0}
\vspace{0.6cm}

Following the preparations from the previous two sections, we are now in a position to complete the
proofs of Theorems 1 to 3. In fact, we obtain more general versions of Theorems 1 and 3, but first we prove Theorem~\ref{satz:sobolev_chain_lp} and inequalities (\ref{bvstrongest}).

\vspace{0.3cm}

\noindent {\it Proof of Theorem~\ref{satz:sobolev_chain_lp} and (\ref{bvstrongest}).}
First, suppose that $1 \leq p < n$ and that $f \in \dot{W}^{1,p}(\mathbb{R}^n)$. We may also assume that $f$ is not identically 0.
Next, note that for every $K \in \mathcal{K}^n$ containing the origin in its interior and $\lambda > 0$, we have $S_p(\lambda K,\cdot) = \lambda^{n-p}S_p(K,\cdot)$.
Consequently, taking
\[K = \left | \langle f \rangle_p \right |^{-\frac{1}{n-p}} \langle f \rangle_p,  \]
it follows from Theorem \ref{the:lyz_volume_normalized} that for $1 \leq i \leq n$ and $E \in \mathrm{Gr}_{n,i}$,
\[q_{i,p}\int_{\mathbb{R}^n} \|\nabla f(x) |E \|^p \,dx = q_{i,p}\int_{\mathbb{S}^{n-1}} \| u|E \|^p \,dS_p(K,u).\]
Hence, definition (\ref{defeip}) and Theorem \ref{the:geometric_inequality} yield that for $1 \leq i < j \leq n$,
\[\mathcal{E}_{j,p} \geq  \mathcal{E}_{i,p} \]
which proves Theorem~\ref{satz:sobolev_chain_lp}.

The inequalities from (\ref{bvstrongest}) for $f \in BV(\mathbb{R}^n)$ not identically 0, follow from similar arguments, by taking $K =  \langle f \rangle$ and applying Theorems \ref{the:lyz_bounded_variation} and \ref{the:geometric_inequality}. Alternatively, (\ref{bvstrongest}) can also be deduced from Theorem \ref{satz:sobolev_chain_lp} by an approximation argument. \hfill $\blacksquare$

\vspace{0.3cm}

Note that a combination of Theorem \ref{satz:sobolev_chain_lp} with the affine $L^p$ Sobolev inequality (\ref{lpaffSob}), directly implies the family of $L^p$ Sobolev inequalities (\ref{eq:sobolev_lp}) from Theorem \ref{satz:sobolev_lp}. Similarly, inequalities (\ref{bvstrongest}) and Wang's extension of the affine $L^1$ Sobolev inequality to $BV(\mathbb{R}^n)$ yield the family of Sobolev inequalities (\ref{bvchainsob}) from Theorem \ref{satz:sobolev}. However, in both cases the equality conditions remain to be settled.

In the following, we are not going to approach the characterization of extremal functions in Theorems 1 to 3 directly, but rather establish generalizations of these theorems that are motivated by the fact that in the classical $L^p$ Sobolev inequality (\ref{lpsob}) the Euclidean norm of the gradient can be replaced by an arbitrary norm on $\mathbb{R}^n$ (see, e.g., \textbf{\cite{cordnazvill}}). To this end, let $1 \leq i \leq n - 1$, $1 \leq p < n$ and suppose that $\mu$ is an even measure on $\mathbb{S}^{n-1}$ such that $\mathrm{span}\, \mathrm{supp}\, \mu = E_i$. For $f \in \dot{W}^{1,p}(\mathbb{R}^n)$, we now define
\begin{equation} \label{defEipmu}
\mathcal{E}_{i,p}^{\mu}(f) = \left ( \int_{\mathrm{Gr}_{n,i}}\left ( \int_{\mathbb{R}^n} \| \nabla f(x)|E\|^p_{Z_p^{\mu}(E)^{\circ}} \,dx \right )^{-n/p}dE\right )^{-1/n},
\end{equation}
where $Z_p^{\mu}$ denotes again the $L^p$ zonoid generated by $\mu$ (see Section 2).

The problem of whether a version of Theorem \ref{satz:sobolev_chain_lp} also holds for more general norms than the Euclidean one was first raised by Ludwig. With our next result we answer this question in the affirmative, when the unit ball of the norm is a polar $L^p$ zonoid.

\begin{theorem} \label{lpsobgeneralnorm}
Suppose that $1\leq i \leq n-1$, $1 \leq p < n$, and that $\mu$ is an even measure on $\mathbb{S}^{n-1}$ such that $\mathrm{span}\, \mathrm{supp}\, \mu = E_i$.
If $f \in \dot{W}^{1,p}(\mathbb{R}^n)$, then
\begin{equation} \label{gennorminequ}
\mathcal{E}_{i,p}^{\mu}(f) \geq \mu(\mathbb{S}^{n-1})^{1/p}c_{n,p}\|f\|_{p^*}
\end{equation}
with equality for $p > 1$, if and only if $f(x)$ has the form (\ref{equallpsob}) when $i > 1$, and if and only if $f(x)$ has the form (\ref{equ1742}) when $i = 1$. Moreover, if $\mu(\mathbb{S}^{n-1}) = 1$, then
\begin{equation} \label{gennormstrong}
\mathcal{E}_{n,p}(f) \geq \mathcal{E}_{i,p}^{\mu}(f) \geq \mathcal{E}_{1,p}(f).
\end{equation}
\end{theorem}

\noindent {\it Proof.} We begin with the proof of (\ref{gennormstrong}). To this end, we may assume that $f$ is not identically 0. Since $1 \leq i \leq n - 1$, we have, by (\ref{suppprojection}) and (\ref{mulifttograss}),
\begin{equation} \label{bigproof17}
\int_{\mathrm{Gr}_{n,i}}\!\!\left ( \int_{\mathbb{R}^n}\!\! \| \nabla f(x)|E\|^p_{Z_p^{\mu}(E)^{\circ}} \,dx\! \right )^{\!\!-n/p}\!dE =
\int_{\mathrm{SO}(n)}\!\!\left ( \int_{\mathbb{R}^n}\!\! h(\phi Z_p^{\mu},\nabla f(x))^p \,dx\! \right )^{\!\!-n/p}\!d\phi.
\end{equation}
Hence, by definition (\ref{defEipmu}) and Jensen's inequality,
\[\mathcal{E}_{i,p}^{\mu}(f) \leq \left (  \int_{\mathrm{SO}(n)}\!\int_{\mathbb{R}^n}\!\! h(\phi Z_p^{\mu},\nabla f(x))^p \,dx\,d\phi \right )^{\!1/p}.  \]
Thus, by (\ref{eq:support_function_thetainv}) and (\ref{eq:rotation_one_dim}), Fubini's theorem, the invariance of the Haar measure on $\mathrm{SO}(n)$, and the fact that $\mu(\mathbb{S}^{n-1}) = \breve{\mu}(\mathrm{SO}(n)) = 1$, we have
\begin{align*}
\mathcal{E}_{i,p}^{\mu}(f) & \leq \left ( \int_{\mathbb{R}^n}\! \int_{\mathrm{SO}(n)}\! \int_{\mathrm{SO}(n)}\!\!\! |\nabla f(x) \cdot \phi \psi e_n|^p \,d\phi\,d\breve{\mu}(\psi)\,dx \right )^{\!1/p} \\
& = \left ( \int_{\mathbb{R}^n}\! \int_{\mathrm{SO}(n)}\!\!\! |\nabla f(x) \cdot \phi e_n|^p \,d\phi\,dx \right )^{\!1/p}
= \left (\! q_{n,p}\!\int_{\mathbb{R}^n}\!\! \|\nabla f(x)\|^p \,dx \right )^{\!1/p} = \mathcal{E}_{n,p}(f),
\end{align*}
which proves the left-hand inequality in (\ref{gennormstrong}).

In order to prove the right-hand inequality in (\ref{gennormstrong}), we use as before (\ref{bigproof17}), (\ref{eq:support_function_thetainv}), and (\ref{eq:rotation_one_dim}), followed by Fubini's theorem, to see that
\begin{equation} \label{eipmudarst}
\mathcal{E}_{i,p}^{\mu}(f)= \left (\int_{\mathrm{SO}(n)}\!\! \left ( \int_{\mathrm{SO}(n)}\! \int_{\mathbb{R}^n}\! |\nabla f(x) \cdot \phi \psi e_n|^p \,dx\,d\breve{\mu}(\psi)\!\right )^{\!-n/p}\!d\phi \!\right )^{\!-1/n}.
\end{equation}
Since $\breve{\mu}(\mathrm{SO}(n)) = 1$, an application of Jensen's inequality therefore shows that
\begin{equation} \label{equality1717}
\mathcal{E}_{i,p}^{\mu}(f) \geq \left (\int_{\mathrm{SO}(n)}\!  \int_{\mathrm{SO}(n)}\!\! \left ( \int_{\mathbb{R}^n}\! |\nabla f(x) \cdot \phi \psi e_n|^p \,dx \!\right )^{\!-n/p} \!d\breve{\mu}(\psi)\,d\phi \!\right )^{\!-1/n}.
\end{equation}

\pagebreak

\noindent Finally, Fubini's theorem, the invariance of the Haar measure on $\mathrm{SO}(n)$, (\ref{mulifttograss}), and the fact that $q_{1,p}=1$, yield the desired inequality
\begin{align*}
\mathcal{E}_{i,p}^{\mu}(f) & \geq \left (\int_{\mathrm{SO}(n)}\!  \int_{\mathrm{SO}(n)}\!\! \left ( \int_{\mathbb{R}^n}\! |\nabla f(x) \cdot \phi \psi e_n|^p \,dx \!\right )^{\!-n/p} \!d\phi\,d\breve{\mu}(\psi) \!\right )^{\!-1/n} \\
& = \left (\int_{\mathrm{SO}(n)}\!\! \left ( \int_{\mathbb{R}^n}\! \|\nabla f(x)| \phi E_1\|^p \,dx \!\right )^{\!-n/p} \!d\phi \!\right )^{\!-1/n} = \mathcal{E}_{1,p}(f).
\end{align*}

Next, note that the homogeneity of $\mathcal{E}_{i,p}^{\mu}$ and a combination of the right-hand inequality in (\ref{gennormstrong}) with the affine $L^p$ Sobolev inequality (\ref{lpaffSob}), yield inequality (\ref{gennorminequ}). Moreover, since $h(Z_p^{\mu},x)^p$ is a positive multiple of $\|x|E_1\|^p$ when $i = 1$, we see that (\ref{gennorminequ}) reduces to (\ref{lpaffSob}) in this case. In particular, equality holds in (\ref{gennorminequ}) when $i = 1$ if and only if $f(x)$ has the form (\ref{equ1742}).

It remains to settle the equality conditions for (\ref{gennorminequ}) when $2 \leq i \leq n - 1$. To this end, note that by (\ref{gennormstrong}), any function of the form (\ref{equallpsob}) must be an extremizer of (\ref{gennorminequ}). In order to show the converse, we may assume again that $f$ is not identically 0 and that $\mu(\mathbb{S}^{n-1})=1$. First, we define (as in the proof of Theorem \ref{satz:sobolev_chain_lp} above)
\[K = \left | \langle f \rangle_p \right |^{-\frac{1}{n-p}} \langle f \rangle_p.  \]
Then, by Theorem \ref{the:lyz_volume_normalized}, (\ref{eipmudarst}), and the fact that $\mu$ is the pushforward of $\breve{\mu}$ under the projection $\pi$ (see Section 2), we see that
\[\mathcal{E}_{i,p}^{\mu}(f)= \left (\int_{\mathrm{SO}(n)}\!\! \left ( \int_{\mathbb{S}^{n-1}}\! \int_{\mathbb{S}^{n-1}}\!\! |u \cdot \phi v|^p \,dS_p(K,u)\,d\mu(v)\!\right )^{\!-n/p}\!d\phi \!\right )^{\!-1/n}.\]
Now, since, by (\ref{gennormstrong}), equality in (\ref{gennorminequ}) implies equality in (\ref{equality1717}), we must have
\[ \left ( \int_{\mathbb{S}^{n-1}}\! \int_{\mathbb{S}^{n-1}}\!\! |u \cdot \phi v|^p \,dS_p(K,u)\,d\mu(v)\!\right )^{\!-n/p} =
\int_{\mathbb{S}^{n-1}}\! \left ( \int_{\mathbb{S}^{n-1}}\!\! |u \cdot \phi v|^p \,dS_p(K,u)\!\right )^{\!-n/p} \!d\mu(v)  \]
or, equivalently, by (\ref{deflpprojbod}),
\begin{equation} \label{lpprojequ}
\left ( \int_{\mathbb{S}^{n-1}}\! h(\Pi_pK,\phi v)^p\,d\mu(v)\!\right )^{\!-n/p} =
\int_{\mathbb{S}^{n-1}}\! h(\Pi_pK,\phi v)^{-n} \,d\mu(v)
\end{equation}
for every $\phi \in \mathrm{SO}(n)$. We claim that this implies that $\Pi_pK$ is a ball. In order to see this, note that by the equality conditions of Jensen's inequality, (\ref{lpprojequ}) holds if and only if for each $\phi \in \mathrm{SO}(n)$, there exists a $c_{\phi}>0$ such that
\[h(\Pi_pK,\phi v) = c_{\phi} \mbox{ for $\mu$-a.e.\ } v \in \mathbb{S}^{n-1}.  \]
Thus, since $2 \leq i \leq n - 1$ and $\mathrm{span}\, \mathrm{supp}\, \mu = E_i$, there exist (at least) two linearly independent unit vectors $u_1, u_2 \in \mathbb{S}^{n-1}$ such that
for each $\phi \in \mathrm{SO}(n)$,
\begin{equation} \label{babyyoda}
h(\Pi_pK,\phi u_1) = h(\Pi_pK,\phi u_2) = c_{\phi}.
\end{equation}
Let $t = u_1 \cdot u_2$ and denote by $H_{u_1,t} = \{x \in \mathbb{R}^n: u_1 \cdot x = t\}$. Then, $-1 < t <1$ and for $w \in \mathbb{S}^{n-1} \cap H_{u_1,t}$, there exists $\vartheta \in \mathrm{SO}(n)$ such that $\vartheta u_1 = u_1$ and $\vartheta u_2 = w$. Replacing $\phi$ by $\phi \vartheta$ in (\ref{babyyoda}), thus yields
\[c_{\phi\vartheta} = h(\Pi_pK,\phi\vartheta u_1) = h(\Pi_pK,\phi\vartheta u_2) =  h(\Pi_pK,\phi u_1) =  h(\Pi_pK,\phi w) = c_{\phi}.  \]
Since $w \in \mathbb{S}^{n-1} \cap H_{u_1,t}$ was arbitrary, we see that for each $\phi \in \mathrm{SO}(n)$, there exists $c_{\phi}>0$ such that
\[h(\Pi_pK,\phi v) = c_{\phi} \mbox{ for all } v \in \mathbb{S}^{n-1} \cap H_{u_1,t}\]
or, equivalently,
\begin{equation} \label{criticalyoda}
h(\Pi_pK,u) = c_{\phi} \mbox{ for all } u \in \mathbb{S}^{n-1} \cap H_{\phi u_1,t}.
\end{equation}
In particular, by choosing $\phi$ to be the identity, we obtain
\[h(\Pi_p K,u) = c_{\mathrm{id}} \mbox{ for all } u \in \mathbb{S}^{n-1} \cap H_{u_1,t}.\]
Now, if we choose $\phi$ in (\ref{criticalyoda}) such that $\mathbb{S}^{n-1} \cap H_{u_1,t}$ and $\mathbb{S}^{n-1} \cap H_{\phi u_1,t}$ have non-empty intersection, then it follows that $c_{\mathrm{id}} = c_{\phi}$. But, since we can reach any point on $\mathbb{S}^{n-1}$ by finitely many iterations of this procedure, we obtain $h(\Pi_p K,u) = c_{\mathrm{id}}$ for all $u \in \mathbb{S}^{n-1}$, that is, $\Pi_p K$ is a ball as desired.

In conclusion, if equality holds in (\ref{gennorminequ}), then $\Pi_p K$ must be a Euclidean ball and, by the right-hand inequality in (\ref{gennormstrong}), equality must also hold in (\ref{lpaffSob}), that is, $f(x)$ is of the form (\ref{equ1742}). The latter implies that there exists an origin-symmetric ellipsoid $E \subseteq \mathbb{R}^n$ and $x_1 \in \mathbb{R}^n$ such that $f(x) = f^E(x + x_1)$ for a.e.\ $x \in \mathbb{R}^n$ (cf.\ proof of Corollary~4.1 in \textbf{\cite{nguyen}}). Hence, by Proposition \ref{prop:lyz_results} (b), $K$ is a dilate of $E = A\mathbb{B}^n$ for suitable $A \in \mathrm{GL}(n)$. However, since
\[\Pi_p(A\mathbb{B}^n) = |\mathrm{det}\, A|^{1/p}\,A^{-\mathrm{T}} \Pi_p \mathbb{B}^n = |\mathrm{det}\, A|^{1/p}\,A^{-\mathrm{T}} \mathbb{B}^n  \]
for any $A \in \mathrm{GL}(n)$ (cf.\ \textbf{\cite{lutwak2000lp}}) and $\Pi_p K$ is a ball, $E$ must be a Euclidean ball as well. This implies that $f$ is an extremizer of (\ref{lpsob}) (see, e.g., \textbf{\cite{wang2013affine}}). \hfill $\blacksquare$

\vspace{0.3cm}

Let us first note that the classical $L^p$ Sobolev inequality (\ref{lpsob}) together with the special case of Theorem \ref{lpsobgeneralnorm}, where $\mu$ is taken to be the $(i-1)$-dimensional Hausdorff measure on $S^{i-1}$ normalized such that $\mu(S^{i-1}) = q_{i,p}^{-1}$, completes the proof of Theorem~\ref{satz:sobolev_lp}. Next, we want to emphasize that by (\ref{gennormstrong}), each of the $L^p$ Sobolev inequalities from (\ref{gennorminequ}) is stronger than the classical inequality (\ref{lpsob}) and that, in turn, the strongest one among them is the affine $L^p$ Sobolev inequality (\ref{lpaffSob}). Finally, let us remark that it is an open problem whether an inequality like (\ref{gennorminequ}) still holds, when the $L^p$ zonoids $Z_p^{\mu}$ in the definition of $\mathcal{E}_{i,p}^{\mu}$ are replaced by more general convex bodies.

\vspace{0.2cm}

With our final result, we establish an extension of Theorem \ref{lpsobgeneralnorm} for the case $p = 1$ to functions of bounded variations which also generalizes Theorem 3:

\begin{theorem} \label{sobbvgeneralnorm}
Suppose that $1\leq i \leq n-1$ and that $\mu$ is an even measure on $\mathbb{S}^{n-1}$ such that $\mathrm{span}\, \mathrm{supp}\, \mu = E_i$.
If $f \in BV(\mathbb{R}^n)$, then
\begin{equation} \label{finalsobbv17}
\mathcal{E}_{i}^{\mu}(f):= \left ( \int_{\mathrm{Gr}_{n,i}}\!\!\left ( \int_{\mathbb{R}^n}\!\! \| \sigma_f|E\|^p_{Z^{\mu}(E)^{\circ}} \,d|Df|\! \right )^{-n}\!dE\!\right )^{\!\!-1/n}
\geq  \frac{2\omega_{n-1}\mu(\mathbb{S}^{n-1})}{\omega_n^{1-1/n} }\, \|f\|_{\frac{n}{n-1}}
\end{equation}
with equality if and only if $f$ is the multiple of a characteristic function of a ball when $i > 1$ and that of an ellipsoid when $i = 1$. Moreover, if $\mu(\mathbb{S}^{n-1}) = 1$, then
\begin{equation} \label{finalbv17}
\mathcal{E}_n(f) \geq \mathcal{E}_{i}^{\mu}(f) \geq \mathcal{E}_1(f).
\end{equation}
\end{theorem}

\noindent {\it Proof.} Since the proof of (\ref{finalbv17}) is almost verbatim the same as that of (\ref{gennormstrong}) (basically, by replacing $\nabla f$ by $\sigma_f$), we will not repeat it here. Having established (\ref{finalbv17}), the Sobolev inequalities (\ref{finalsobbv17}) follow from Wang's extension of the affine Zhang--Sobolev inequality. Moreover, when $i = 1$, inequality (\ref{finalsobbv17}) reduces to the affine Zhang--Sobolev inequality on $BV(\mathbb{R}^n)$.

In order to settle the equality conditions for (\ref{finalsobbv17}) when $2 \leq i \leq n - 1$, one can show as in the proof of Theorem \ref{lpsobgeneralnorm}, that equality in (\ref{finalsobbv17}) for a function $f \in BV(\mathbb{R}^n)$ not identically 0, implies that $\Pi \langle f \rangle$ is a ball. Moreover, by (\ref{finalbv17}) and the (extended) affine Zhang--Sobolev inequality, $f$ must be a multiple of the characteristic function of an ellipsoid $E = A\mathbb{B}^n$ for suitable $A\in\mathrm{GL}(n)$. Since, by Proposition \ref{prop:lyz_results} (a), we have $\langle \mathds{1}_E \rangle = E$, we infer that $\Pi E$ must be a ball. But, since
\[\Pi (A \mathbb{B}^n) = |\det A |A^{-\mathrm{T}}\Pi \mathbb{B}^n= \omega_{n-1}|\det A |A^{-\mathrm{T}} \mathbb{B}^n
	\end{equation*}
for any $A \in \mathrm{GL}(n)$ (cf.\ \textbf{\cite[\textnormal{Theorem 4.1.5}]{gardner2ed}}), $f$ must actually be a multiple of the characteristic function of a ball. \hfill $\blacksquare$

\vspace{1cm}

\noindent {{\bf Acknowledgments} The authors were supported by the European Research Council (ERC), Project number: 306445, and the Austrian Science Fund (FWF),
Project numbers: Y603-N26 and P31448-N35.}

\begin{small}

\[ \begin{array}{ll} \mbox{Philipp Kniefacz} & \mbox{Franz E. Schuster} \\
\mbox{Vienna University of Technology \phantom{wwwwWW}} & \mbox{Vienna University of Technology} \\ \mbox{philipp.kniefacz@tuwien.ac.at} & \mbox{franz.schuster@tuwien.ac.at}
\end{array}\]

\end{small}

\end{document}